\documentclass[a4paper,11pt]{article}
\usepackage[utf8]{inputenc}
\usepackage[english]{babel}
\usepackage{amsfonts,amssymb,amsthm,amsmath}
\usepackage{url}
\usepackage{enumitem}

\usepackage{graphicx}
\usepackage{epstopdf}
\usepackage{varwidth}

\usepackage{a4wide}
\usepackage{xcolor}
\usepackage[colorlinks=true,linktocpage,pdfpagelabels,
bookmarksnumbered,bookmarksopen]{hyperref}
\definecolor{ForestGreen}{rgb}{0.1,0.6,0.05}
\definecolor{EgyptBlue}{rgb}{0.063,0.1,0.6}
\hypersetup{
	colorlinks=true,
	linkcolor=EgyptBlue,         
	citecolor=ForestGreen,
	urlcolor=olive
}

\usepackage[hyperpageref]{backref}

\let\OLDthebibliography\thebibliography
\renewcommand\thebibliography[1]{
	\OLDthebibliography{#1}
	\setlength{\parskip}{1pt}
	\setlength{\itemsep}{1pt plus 0.3ex}
}

\numberwithin{equation}{section}
\newtheorem{theorem}{Theorem}[section]
\newtheorem{lemma}[theorem]{Lemma}
\newtheorem{prop}[theorem]{Proposition}
\newtheorem{cor}[theorem]{Corollary}
\theoremstyle{definition}
\newtheorem{remark}[theorem]{Remark}

\title{
	\vspace*{-1cm}
On exact Pleijel's constant for some domains
}

\makeatletter
\newcommand{\address}[1]{\gdef\@address{#1}}
\newcommand{\email}[1]{\gdef\@email{\url{#1}}}
\newcommand{\homepage}[1]{\gdef\@homepage{\url{#1}}}
\newcommand{\@endstuff}{\vspace{\baselineskip}\noindent\small
	\begin{tabular}{@{}l}
		\@address\\
		\textit{E-mail address:} \@email\\
		\textit{Homepage:} \@homepage
	\end{tabular}}
\AtEndDocument{\@endstuff}
\makeatother

\author{Vladimir Bobkov}
\email{bobkov@kma.zcu.cz}
\homepage{http://vladimir-bobkov.ru}
\address{
	\textsc{Department of Mathematics and NTIS, Faculty of Applied Sciences},\\ 
	\textsc{University of West Bohemia, Univerzitn\'i 8, 306 14 Plze\v{n}, Czech Republic};\\ 
	\textsc{Institute of Mathematics with Computing Centre of RAS},\\ 
	\textsc{Chernyshevsky str. 112, 450008 Ufa, Russia}}

\date{}

\begin{document}
\maketitle 

\begin{abstract}
	We provide an explicit expression for the Pleijel constant for the planar disk and some of its sectors, as well as for $N$-dimensional rectangles. 
	In particular, the Pleijel constant for the disk is equal to $0.4613019\ldots$
	Also, we characterize the Pleijel constant for some rings and annular sectors in terms of asymptotic behavior of zeros of certain cross-products of Bessel functions.
	
\end{abstract}

\section{Introduction and main results}

Let $\Omega \subset \mathbb{R}^N$ be a bounded domain with the boundary $\partial \Omega$, $N \geq 2$.
Consider the eigenvalue problem 
\begin{equation}\label{eq:D}
\left\{
\begin{aligned}
-\Delta u &= \lambda u &&{\rm in}\ \Omega, \\
u&=0 &&{\rm on }\ \partial \Omega,
\end{aligned}
\right.
\end{equation}
and denote by $\{\lambda_n\}_{n \in \mathbb{N}}$ the sequence of the corresponding eigenvalues,\footnote{Here and below we always assume that $\mathbb{N}$ does not contain $0$.}
$$
0 < \lambda_1 < \lambda_2 \leq \dots \leq \lambda_n \to \infty
\quad
\text{as }
n \to \infty.
$$
For any eigenfunction $\varphi_n$ associated with $\lambda_n$, let $\mu(\varphi_n)$ be a number of nodal domains of $\varphi_n$, that is, a number of connected components of $\Omega \setminus \overline{\{x\in\Omega: \varphi_n(x)=0\}}$.
Courant's nodal domain theorem \cite{courant} asserts that $\mu(\varphi_n) \leq n$ for any $n \in \mathbb{N}$.
This result was refined by Pleijel for $N=2$ and then extended by B\'erard \& Meyer to the general $N \geq 2$ as follows.
\begin{theorem}[{\cite[Section 5]{pleijel}} and  {\cite[Theorem II.7]{BM}}]\label{thm:Pleijel}
	Let $j_{\frac{N}{2}-1,1}$ be the first zero of the Bessel function $J_{\frac{N}{2}-1}$. 
	Then
	$$
	Pl(\Omega) := \limsup_{n \to \infty} \frac{\mu(\varphi_n)}{n} \leq \gamma(N) := \frac{2^{N-2} \, N^2 \, \Gamma\left(\frac{N}{2}\right)^2}{j_{\frac{N}{2}-1,1}^N} < 1.
	$$
	In particular, if $N=2$, then 
	\begin{equation}\label{eq:Pleijel}
	Pl(\Omega) \leq \frac{4}{j_{0,1}^2} = 0.6916602\ldots
	\end{equation}
\end{theorem}
We will call $Pl(\Omega)$ \textit{the Pleijel constant} for $\Omega$.
Note that the function $N \mapsto \gamma(N)$ is strictly decreasing and $\lim\limits_{N \to \infty} \frac{\gamma(N+1)}{\gamma(N)} = \frac{2}{e}$, see Helffer \& Sundqvist \cite[Theorem 5.1 and Remark 5.4]{HS}. 

The upper bound $\frac{4}{j_{0,1}^2}$ in \eqref{eq:Pleijel} is not sharp, as it was pointed out by Polterovich \cite{pol} and then rigorously developed by Bougain \cite{B} and Steinerberger \cite{stein}.
Furthermore, Polterovich conjectured that the optimal upper bound for $N=2$ should be
$$
Pl(\Omega) \leq \frac{2}{\pi} = 0.6366197\ldots
$$
If this upper estimate holds true, then its optimality follows from the consideration of any rectangle $\mathcal{R}(a, b) = (0,a) \times (0,b)$ with irrational ratio $\frac{a^2}{b^2}$, see Helffer \& Hoffmann-Ostenhof \cite[Proposition 5.1]{Helff}. 

We start by extending \cite[Proposition 5.1]{Helff} to the general $N$-dimensional case as follows.
\begin{prop}\label{prop:rect}
	Let $\mathcal{R}(a_1,\dots,a_N) = (0,a_1)\times\cdots\times(0,a_N)$ be a $N$-orthotope such that $\frac{a_i^2}{a_j^2}$ is irrational for any $i \neq j$.
	Then
	$$
	Pl(\mathcal{R}(a_1,\dots,a_N)) = \rho(N) := \frac{2^N \Gamma\left(\frac{N}{2}+1\right)}{\pi^\frac{N}{2} N^\frac{N}{2}}.
	$$
	Moreover, the function $N \mapsto \rho(N)$ is strictly decreasing and $\lim\limits_{N \to \infty} \frac{\rho(N+1)}{\rho(N)} = \sqrt{\frac{2}{\pi e}}$.
\end{prop}

It is clear that $\rho(N) < \gamma(N)$ for any $N \geq 2$ by the same reasoning as in \cite[Remark 2.5]{pol}. Furthermore, by analogy with the conjecture of Polterovich \cite{pol}, one can ask whether the estimate $Pl(\Omega) \leq \rho(N)$ is valid for any $\Omega \subset \mathbb{R}^N$ and $N \geq 2$.

\smallskip

Note that the ``irrational'' rectangles discussed above are the most simple domains whose Pleijel's constant can be worked out explicitly. 
However, to the best of our knowledge, the exact value of the Pleijel constant has not been known for any other domains,\footnote{
	For the quantum harmonic oscillator $-\Delta + \sum_{i=1}^N a_i^2 x_i^2$ with rationally independent $a_i$'s, Charron \cite{charron} proved that the Pleijel constant is equal to $N!/N^N$}
and the question of finding of such domains was proposed by Bonnaillie-No\"el et al \cite[Section 6.1]{BHHO}.
Furthermore, numerical experiments even for the planar disk and square (see, e.g., Blum et al \cite{BGS}) demonstrate relatively slow convergence of ratio $\frac{\mu(\varphi_n)}{n}$ as $n \to \infty$. 
Thus, theoretical results in this direction are of some interest.

The aim of the present note is to obtain the explicit expression for the Pleijel constant for the planar disk. 
Our main result is the following.

\begin{theorem}\label{thm:disk}
	Let $B := \{x \in \mathbb{R}^2:~ |x|<1\}$. Then 
	\begin{equation*}\label{eq:PLB0}
	Pl(B) = 8 \, \sup_{x>0} \left\{ x \left(\cos \theta(x)\right)^2 \right\}  = 0.4613019\ldots,
	\end{equation*}
	where $\theta=\theta(x)$ is the solution of the transcendental equation
	\begin{equation*}\label{eq:theta0}
	\tan \theta - \theta = \pi x, 
	\quad 
	\theta \in \left(0, \frac{\pi}{2}\right).
	\end{equation*}
\end{theorem}

Note that the value of $Pl(B)$ is in the good correspondence with the numerical simulation of Blum et al \cite[Fig.~1]{BGS}. 
We also refer the reader to Han et al \cite{HMT} where the sharp estimates for the length of the nodal set of eigenfunctions of \eqref{eq:D} on $B$ are obtained.

Consider now the circular sector with the angle $\alpha \in (0,2\pi]$ defined as
$$
\Sigma^\alpha := 
\left\{(\varrho \cos\vartheta, \varrho \sin\vartheta):~ \varrho \in (0,1),~ \vartheta \in \left(-\frac{\alpha}{2},\frac{\alpha}{2}\right)
\right\}.
$$

\begin{theorem}\label{thm:sector-disk}
	Let $\alpha$ be such that any eigenvalue 
	of \eqref{eq:D} on $\Sigma^\alpha$ has the multiplicity $1$.
	Then $Pl(\Sigma^\alpha) = Pl(B)$.
	In particular, the assumption on $\alpha$ is satisfied if $\alpha = \pi/m$, $m\in \mathbb{N}$.
\end{theorem}

\begin{remark}
	There exists $\alpha \in (0, 2\pi)$ such that there is an eigenvalue of \eqref{eq:D} on $\Sigma^\alpha$ whose multiplicity is at least $2$, see Bonnaillie-No\"el \& L\'ena \cite[Section 2.3]{BL}. 
	In general, the question about the multiplicity $1$ for eigenvalues of \eqref{eq:D} on $\Sigma^\alpha$ is equivalent to the following: 
	Find assumptions on $\alpha \in (0,2\pi]$ which guarantee that the Bessel functions $J_{\nu_1 \frac{\pi}{\alpha}}$ and $J_{\nu_2 \frac{\pi}{\alpha}}$ have no common positive zeros for any  $\nu_1, \nu_2 \in \mathbb{N}$.	
	This question is reminiscent of Bourget's hypothesis \cite[p.~484]{watson}.
\end{remark}

Finally, we characterize the Pleijel constant for rings and annular sectors. 
First, consider the ring (annulus) $A_{r} := \{x \in \mathbb{R}^2:~ r<|x|<1\}$, $r \in (0,1)$.
\begin{prop}\label{thm:annulus}
	Let $r$ be such that any eigenvalue of \eqref{eq:D} on $A_{r}$ has the multiplicity at most $2$.
	Then
	\begin{equation*}\label{eq:PLBa}
	Pl(A_{r}) = \frac{8}{1-r^2} \, \sup_{x>0} \left\{x \, \limsup_{k \to \infty} \frac{k^2}{a_{kx,k}^2}\right\},
	\end{equation*}
	where $a_{kx,k}$ is the $k$-th positive zero of the cross-product of Bessel functions
	$J_{kx}(r z) Y_{kx}(z) - J_{kx}(z) Y_{kx}(r z)$.
\end{prop}

Consider now the annular sector with the inner radius $r \in (0,1)$ and the angle $\alpha \in (0,2\pi]$ defined as
$$
\Sigma^\alpha_{r} := 
\left\{(\varrho \cos\vartheta, \varrho \sin\vartheta):~ \varrho \in (r,1),~ \vartheta \in \left(-\frac{\alpha}{2},\frac{\alpha}{2}\right)
\right\}.
$$
\begin{prop}\label{thm:annulus_sector}
	Let $r$ and $\alpha$ be such that any eigenvalue of \eqref{eq:D} on $\Sigma^\alpha_{r}$ has the multiplicity $1$.
	Then $Pl(\Sigma^\alpha_{r}) = Pl(A_r)$.
\end{prop}

\begin{remark}
	It seems that, by now, the asymptotic behavior of $a_{kx,x}$ as $k \to \infty$ is not studied as comprehensively as whose of $j_{kx,k}$ (\cite{EL1,E}), which obstructs us to obtain an explicit expression for $Pl(A_{r})$ in the spirit of Theorem \ref{thm:disk}. 
	We expect that the following asymptotics should be valid:
	$$
	Pl(A_{r}) \to Pl(B)
	\quad \text{as } r \to 0;
	\qquad 
	Pl(A_{r}) \to \rho(2) = \frac{2}{\pi}
	\quad \text{as } r \to 1.
	$$
\end{remark}

Let us note that there exist parameters $r$ and $\alpha$ for which the assumptions of Propositions \ref{thm:annulus} and \ref{thm:annulus_sector} do not hold.
\begin{lemma}\label{lem:multiplicity}
	There exists $r$ (and $\alpha$) such that there is an eigenvalue of \eqref{eq:D} on $A_{r}$ (on $\Sigma^\alpha_{r}$) whose multiplicity is at least $3$ (at least $2$), see Fig.~\ref{fig1}.
\end{lemma}
But it should be also emphasized that we do not know whether the sets of admissible parameters for Propositions \ref{thm:annulus} and \ref{thm:annulus_sector} are non-empty.


\section{Proofs}

The structure of all the proofs consists of three main steps. First, we characterize eigenvalues and eigenfunctions of \eqref{eq:D} explicitly as multi-indexed sets, due to the ``separable'' nature of the considered domains. Then, after reordering the set of eigenvalues into the increasing sequence $\{\lambda_n\}_{n \in \mathbb{N}}$, we apply the Weyl law (see, e.g., Ivrii \cite{ivrii}) to express the labeling $n$ in terms of the corresponding eigenvalue $\lambda_n$ for sufficiently large $n$. Finally, we elaborate the obtained expression for $Pl(\Omega)$ up to the desired form.

\subsection{Proof of Proposition \ref{prop:rect}}
Let $N \geq 2$ and let $\mathcal{R}(a_1,\dots,a_N) = (0,a_1)\times\cdots\times(0,a_N)$ be a $N$-orthotope such that the ratio $\frac{a_{i}^2}{a_j^2}$ is irrational for any $i \neq j$.
By separation of variables, any eigenvalue of \eqref{eq:D} on $\mathcal{R}(a_1,\dots,a_N)$ is given by 
$$
\lambda_{m_1,\dots,m_N} = \frac{\pi^2 m_1^2}{a_1^2} + \dots + \frac{\pi^2 m_N^2}{a_N^2}, \quad
m_1,\dots,m_N \in \mathbb{N},
$$
with the associated eigenfunction
$$
\varphi_{m_1,\dots,m_N} = \sin \frac{\pi m_1 x}{a_1} \cdots \sin \frac{\pi m_N x}{a_N}.
$$
Evidently, $\mu(\varphi_{m_1,\dots,m_N}) = m_1 \cdots m_N$.

Since $\frac{a_{i}^2}{a_j^2}$ is irrational for any $i \neq j$, all eigenvalues are simple (have the multiplicity $1$). Therefore, no other eigenfunctions occur, and if we put $\lambda_{m_1,\dots,m_N}$'s in the increasing order as $\lambda_1 \leq \dots \lambda_n \leq \dots$, then for any $n$ there exists a unique $N$-tuple $(m_{1,n},\dots,m_{N,n})$ such that $\lambda_n = \lambda_{m_{1,n},\dots,m_{N,n}}$.
Since the explicit relation between $n$ and $(m_{1,n},\dots,m_{N,n})$ is not clear, we use the Weyl law which states that
$$
n = (2\pi)^{-N} \, \omega_N\, |\mathcal{R}(a_1,\dots,a_N)| \, \lambda_n^{N/2} + o(\lambda_n^{N/2}),
$$
where $\omega_N$ is the volume of the unit ball in $\mathbb{R}^N$ and $|\mathcal{R}(a_1,\dots,a_N)|$ is the volume of $\mathcal{R}(a_1,\dots,a_N)$. 
Thus, we get
$$
n = \frac{\omega_N\, a_1 \cdots a_N}{2^{N}}  \left(\frac{m_{1,n}^2}{a_1^2} + \dots + \frac{m_{N,n}^2}{a_N^2}\right)^{N/2} + o(\lambda_n^{N/2}).
$$
Substituting this expression into the definition of the Pleijel constant and omitting, for simplicity, the subindex $n$ in $m_{i,n}$, we obtain
$$
Pl(\mathcal{R}(a_1,\dots,a_N)) 
= 
\frac{2^N}{\omega_N} \limsup_{m_1+\cdots+m_N \to \infty} 
\frac{\frac{m_1}{a_1} \cdots \frac{m_N}{a_N}}{\left(\frac{m_{1}^2}{a_1^2} + \dots + \frac{m_{N}^2}{a_N^2}\right)^{N/2}}.
$$
Applying the inequality between arithmetic and geometric means, we deduce that 
$$
Pl(\mathcal{R}(a_1,\dots,a_N)) 
\leq 
\frac{2^N}{\omega_N N^\frac{N}{2}}.
$$
On the other hand, taking a sequence of $N$-tuples $\{(m_{1,k},\dots,m_{N,k})\}_{k \in \mathbb{N}} \subset \mathbb{N}^N$ such that $\lim\limits_{k \to \infty} \frac{m_{i,k}}{k} = a_i$ for each $i =1,\dots,N$, we get
$$
Pl(\mathcal{R}(a_1,\dots,a_N)) 
\geq 
\frac{2^N}{\omega_N N^\frac{N}{2}},
$$
and hence the equality holds. Since $\omega_N = \frac{\pi^\frac{N}{2}}{\Gamma\left(\frac{N}{2}+1\right)}$, we conclude that
$$
Pl(\mathcal{R}(a_1,\dots,a_N)) 
= 
\rho(N) := 
\frac{2^N \, \Gamma\left(\frac{N}{2}+1\right)}{\pi^\frac{N}{2} \, N^\frac{N}{2}}.
$$

Let us show now that $\frac{\rho(N+1)}{\rho(N)} < 1$ for any $N \geq 2$. We have
\begin{equation}\label{eq:rr<1}
\frac{\rho(N+1)}{\rho(N)} = \frac{2\, N^\frac{N}{2} \Gamma\left(\frac{N}{2} + \frac{1}{2} +1\right)}{\pi^\frac{1}{2} \, (N+1)^\frac{N+1}{2}\Gamma\left(\frac{N}{2}+1\right)}.
\end{equation}
Applying Gautschi's-type inequality $\Gamma(x+1) < \left(x+\frac{1}{2}\right)^\frac{1}{2}\Gamma\left(x+\frac{1}{2}\right)$ which holds for all $x > 0$ (see Kershaw \cite[(1.3)]{kershaw}), we have
$$
\Gamma\left(\frac{N}{2} + \frac{1}{2} +1\right) < \left(\frac{N+2}{2}\right)^{1/2} 
\Gamma\left(\frac{N}{2}+1\right).
$$
Substituting this inequality into \eqref{eq:rr<1} and estimating $N^\frac{N}{2} < (N+1)^\frac{N}{2}$, we obtain 
$$
\frac{\rho(N+1)}{\rho(N)} < \left(\frac{2 (N+2)}{\pi (N+1)}\right)^{1/2} < 1
$$
for any $N \geq 2$, that is, the function $N \mapsto \rho(N)$ is strictly decreasing. 

Finally, we show that $\lim\limits_{N\to\infty} \frac{\rho(N+1)}{\rho(N)} = \sqrt{\frac{2}{\pi e}}$.
Let us rewrite \eqref{eq:rr<1} as
$$
\frac{\rho(N+1)}{\rho(N)} = \frac{2^\frac{1}{2}}{\pi^\frac{1}{2}} \left(1+ \frac{1}{N}\right)^{-\frac{N}{2}}  \frac{\left(N+2\right)^\frac{1}{2}}{(N+1)^\frac{1}{2}} \frac{\Gamma\left(\frac{N}{2} + \frac{1}{2} +1\right)}{\Gamma\left(\frac{N}{2}+1\right) \left(\frac{N}{2}+1\right)^\frac{1}{2}}.
$$
Then, noting that $\lim\limits_{x\to\infty}\frac{\Gamma(x+1)}{\Gamma\left(x+\frac{1}{2}\right)\left(x+\frac{1}{2}\right)^\frac{1}{2}} = 1$ (see, e.g., \cite[(2.3)]{kershaw}) and $\lim\limits_{N\to\infty}  \left(1+ \frac{1}{N}\right)^{-\frac{N}{2}} = e^{-\frac{1}{2}}$, we get the desired result.
\qed

\subsection{Proof of Theorem \ref{thm:disk}}

Consider the unit planar disk $B = \{x \in \mathbb{R}^2:~ |x| < 1\}$. 
It is well-known that, by separation of variables, there exists a basis of eigenfunctions of \eqref{eq:D} on $B$ expressed (up to rotation) in polar coordinates $(\varrho,\vartheta)$ as 
$$
\varphi_{\nu,k}(\varrho,\vartheta) = J_\nu(j_{\nu,k} \varrho) \cos(\nu \vartheta),
\quad 
\nu \in \mathbb{N} \cup \{0\},~
k \in \mathbb{N},
$$
and $\lambda_{\nu,k} = j_{\nu,k}^2$ is the eigenvalue associated with $\varphi_{\nu,k}$. Here $j_{\nu,k}$ stands for the $k$-th positive zero of the Bessel function $J_\nu$.
Any eigenvalue $\lambda_{0,k}$ has the multiplicity $1$ (and $\varphi_{0,k}$ is radial), while any other eigenvalue has the multiplicity $2$ (and corresponding eigenfunctions are nonradial, one is a rotation of another). These facts follow from the validity of Bourget's hypothesis, which asserts that $J_\nu$ and $J_{\nu+m}$ do not have common positive zeros for any natural $m$, see \cite[p.~484]{watson}.

Clearly, $\mu(\varphi_{0,k}) = k$ and $\mu(\varphi_{\nu,k}) = 2 \nu k$ for $\nu \in \mathbb{N}$. Below, for brevity, we combine both cases by writing $\mu(\varphi_{\nu,k}) = (2 \nu + \sigma(\nu))k$, where $\sigma(0)=1$ and $\sigma(\nu)=0$ for $\nu \in \mathbb{N}$.

Let us put $\lambda_{\nu,k}$'s in the increasing order as $\lambda_1 \leq \dots \leq \lambda_n \leq \dots$. 
Since we are interested in the behavior as $n \to \infty$ and the explicit relation between $(\nu,k)$ and $n$ is not known, we use the Weyl law to get
$$
n = \lambda_n \frac{|B|^2}{4 \pi^2} + o(\lambda_n) = \frac{\lambda_n}{4} + o(\lambda_n).
$$
Hence, noting that for any $\lambda_n$ there exists a unique pair $(\nu_n,k_n)$ such that $\lambda_n = \lambda_{\nu_n,k_n} = j_{\nu_n,k_n}^2$, we obtain
$$
n = \frac{j_{\nu_n,k_n}^2}{4} + o(j_{\nu_n,k_n}^2).
$$
Substituting this relation to the definition of $Pl(B)$, we deduce that
$$
Pl(B) = \limsup_{n \to \infty} \frac{4 (2 \nu_n + \sigma(\nu_n)) k_n}{j_{\nu_n,k_n}^2}.
$$
Extracting a subsequence which delivers the value $Pl(B)$, omitting (for simplicity) the subindex for $(\nu_n,k_n)$, and noting that $n \to \infty$ iff $\nu + k \to \infty$, we obtain
$$
Pl(B) = \lim_{\nu + k \to \infty} \frac{4 (2 \nu + \sigma(\nu)) k}{j_{\nu,k}^2}.
$$

All we need now is to study the behavior of $j_{\nu,k}$ as $\nu + k \to \infty$. We will use the inequality of McCann \cite[Corollary, p.~102]{mccann} which states that
\begin{equation}\label{eq:mccann}
j_{\nu,k} > \left(\nu^2 + \pi^2 \left(k-\frac{1}{4}\right)^2\right)^{1/2}
\quad 
\text{for any } \nu \geq 0
\text{ and } k \in \mathbb{N}.
\end{equation}

Note first that the sequence of $\varphi_{0,k}$'s cannot be a maximizing sequence for $Pl(B)$ since otherwise the inequality \eqref{eq:mccann} implies $j_{0,k} > \pi(k-1)$ and hence
$$
Pl(B)  = \lim_{k\to\infty} \frac{4k}{j_{0,k}^2} \leq \lim_{k\to\infty} \frac{4k}{\pi^2(k-1)^2} = 0,
$$
but we will see later that $Pl(B) > 0$. 
Thus, we may assume that $\nu > 0$, which yields  $\mu(\varphi_{\nu,k}) = 2 \nu k$ and 
\begin{equation}\label{1}
Pl(B) = \lim_{\nu + k \to \infty} \frac{8 \nu k}{j_{\nu,k}^2}.
\end{equation}

Estimating now \eqref{eq:mccann} from below as $j_{\nu,k} > \nu$ and as $j_{\nu,k} > \frac{\pi k}{2}$, we deduce that 
\begin{equation}\label{eq:upper_bound}
Pl(B) \leq \lim_{\nu + k \to \infty} \min\left\{
\frac{8k}{\nu}, \frac{32 \nu}{\pi^2 k}.
\right\}
\end{equation}
As we already mentioned, it will be shown later that $Pl(B) > 0$.
Using this inequality, we conclude from \eqref{eq:upper_bound} that both $\nu$ and $k$ tend to infinity, and there exist $A_1, A_2 > 0$ such that 
$$
A_1 \nu < k < A_2 \nu
\quad 
\text{for all sufficiently large }
\nu \in \mathbb{N}.
$$
Recalling that $(\nu,k)$ is a maximizing subsequence for $Pl(B)$, we can select a sub-subsequence (which is hence also a maximizing subsequence for $Pl(B)$) still denoted by $(\nu,k)$, such that
\begin{equation}\label{2}
\lim_{\nu \to \infty}\frac{k}{\nu} = x_0 \in [A_1,A_2].
\end{equation}
That is, we have $k = \nu \, x_0 + o(\nu)$ for all large $\nu \in \mathbb{N}$.

Let us now use the result of Elbert \& Laforgia \cite{EL1} (see \cite[Section 1.5]{E} for the precise statement employed in \eqref{3}) which states that	
\begin{equation}\label{3}
\lim_{\nu \to \infty} \frac{j_{\nu,\nu x}}{\nu} = \frac{1}{\cos \theta(x)}, 
\quad x > 0,
\end{equation}
where $\theta=\theta(x)$ is the solution of the (transcendental) equation
\begin{equation}\label{eq:theta}
\tan \theta - \theta = \pi x, 
\quad 
\theta \in \left(0, \frac{\pi}{2}\right).
\end{equation}
Combining \eqref{1}, \eqref{2}, and \eqref{3}, we see that $Pl(B) = 8 x_0 \left(\cos \theta(x_0)\right)^2$, and $x_0$ have to satisfy
\begin{equation}\label{eq:PLB}
Pl(B) = 8 x_0 \left(\cos \theta(x_0)\right)^2 = 8 \, \sup_{x>0} \left\{ x \left(\cos \theta(x)\right)^2 \right\} > 0.
\end{equation}

Most likely, \eqref{eq:theta} and hence \eqref{eq:PLB} cannot be resolved in closed forms. However, one can convince himself that the left-hand side of \eqref{eq:theta} is strictly increasing in $\left(0, \frac{\pi}{2}\right)$, and hence the unique root of \eqref{eq:theta} and the value of $Pl(B)$ can be found with arbitrary precision via the standard numerical methods. 
In particular, using the build-in methods of \textsl{Mathematica}, we obtain
$$
Pl(B) = 0.4613019\ldots
\quad 
\text{and}
\quad 
x_0 = \lim_{\nu \to \infty}\frac{k}{\nu} = 0.3710096\ldots
\pushQED{\qed}
\qedhere
$$

\medskip
\begin{remark}
	\textsl{Mathematica}'s code for finding $Pl(B)$ via \eqref{eq:PLB} can look like that:
	\begin{center}
	\begin{varwidth}{\linewidth}
	\begin{verbatim}
	T[x_?NumericQ] := y /. FindRoot[Tan[y] - y == Pi*x, {y, Pi/4}]; 
	FindMaximum[8*x*(Cos[T[x]])^2, {x, 0.37}]
	\end{verbatim}
	\end{varwidth}
	\end{center}
\end{remark}

\smallskip
\begin{remark}
	If we consider the Neumann eigenvalues instead of the Dirichlet ones, then the result of Theorem \ref{thm:disk} remains valid. Indeed, Neumann eigenfunctions are expressed as
	$$
	\phi_{\nu,k}(\varrho,\vartheta) = J_\nu(j'_{\nu,k} \varrho) \cos(\nu \vartheta),
	\quad 
	\nu \in \mathbb{N} \cup \{0\},~
	k \in \mathbb{N},
	$$
	where $j'_{\nu,k}$ is the $k$-th positive zero of the derivative $J_\nu'$ of the Bessel function $J_\nu$. Moreover, $\lambda_{\nu,k} = (j'_{\nu,k})^2$ is the associated eigenvalue. (Note that $\lambda=0$ is an additional eigenvalue corresponding to the constant eigenfunction.)
	According to the result of Ashu \cite[Theorem 3.2]{ashu} (see also a generalization obtained by Helffer \& Sundqvist \cite[Lemma 2.5]{HS}), $J_\nu'$ and $J_{\nu+m}'$ have no common zeros for any $m \in \mathbb{N}$. Thus, the multiplicity of $\lambda_{0,k}$ is $1$, and the multiplicity of any other eigenvalue is $2$. 
	Note that the Weyl law is still valid for the Neumann case (see, e.g., Ivrii \cite{ivrii}).
	Moreover, $j_{0,k} < j_{0,k}' < j_{0,k+1}$ for $k \geq 1$, and $j_{\nu,k-1} < j_{\nu,k}' < j_{\nu,k}$ for $\nu > 0$ and $k \geq 2$. These inequalities	imply that $\phi_{0,k}$ ($\phi_{\nu,k}$ for $\nu>1$) is, in fact, a restriction to $B$ of the Dirichlet eigenfunction $\varphi_{0,k+1}$ ($\varphi_{\nu,k}$ for $\nu>1$) defined on a slightly bigger ball $B_R$.
	Hence, $\mu(\phi_{0,k}) = k+1$ and $\mu(\phi_{\nu,k}) = 2 \nu k$ for $\nu \in \mathbb{N}$. 
	Arguing now as in the proof of Theorem \ref{thm:disk} and estimating $j_{\nu,k}'$ from both sides as above, we obtain the desired result.
	
	Note that the Pleijel theorem (Theorem \ref{thm:Pleijel}) for the Neumann eigenvalues was proved by Polterovich \cite{pol} for $N=2$ and piecewise real analytic domains, and by L\'ena \cite{lena} for $N \geq 2$ and $C^{1,1}$-smooth domains.
\end{remark}

\subsection{Proof of Theorem \ref{thm:sector-disk}}
	Consider the circular sector of a unit disk
	$$
	\Sigma^\alpha = 
	\left\{(\varrho \cos\vartheta, \varrho \sin\vartheta):~ \varrho \in (0,1),~ \vartheta \in \left(-\frac{\alpha}{2},\frac{\alpha}{2}\right)
	\right\},
	\quad
	\alpha \in (0,2\pi].
	$$
	By separation of variables, it can be shown (see, e.g., \cite[Proposition 2.1]{BL}) that there is a basis of eigenfunctions of \eqref{eq:D} on $\Sigma^\alpha$ of the form
	$$
	\varphi_{\nu,k}^\alpha(\varrho,\vartheta) = J_{\nu\frac{\pi}{\alpha}}\left(j_{\nu\frac{\pi}{\alpha},k} \varrho\right) 
	\sin\left(\nu \pi \left(\frac{\vartheta}{\alpha}+\frac{1}{2}\right)\right),
	\quad 
	\nu, k \in \mathbb{N},
	$$
	with the associated eigenvalues $\lambda_{\nu,k} = j_{\nu\frac{\pi}{\alpha},k}^2$. 
	Evidently, $\mu(\varphi_{\nu,k}^\alpha) = \nu k$.
	We assume that $\alpha$ is chosen in such a way that each $\lambda_{\nu,k}$ has the multiplicity $1$, that is, no other eigenfunctions occur.
	
	The proof of the claim follows the same lines as the proof of Theorem \ref{thm:disk}. Let us present the arguments sketchily. Since $|\Sigma^\alpha| = \frac{\alpha}{2}$, the Weyl law implies
	$$
	n = \lambda_n \frac{|\Sigma^\alpha| |B|}{4 \pi^2} + o(\lambda_n) = \frac{\lambda_n \alpha}{8 \pi} + o(\lambda_n)
	=
	\frac{j_{\nu_n \frac{\pi}{\alpha},k_n}^2 \alpha}{8\pi} + o(j_{\nu_n \frac{\pi}{\alpha},k_n}^2),
	$$
	and hence, for a maximizing subsequence, 
	\begin{equation*}\label{1_2}
	Pl(\Sigma^\alpha) 
	= \lim_{\nu+k \to \infty} \frac{8 \pi \nu k}{j_{\nu \frac{\pi}{\alpha},k}^2 \alpha}
	= \lim_{\tau + k \to \infty} \frac{8 \tau k}{j_{\tau,k}^2},
	\end{equation*}
	where $\tau = \nu \frac{\pi}{\alpha}$.
	Applying the inequality \eqref{eq:mccann}, we deduce that, under the assumption $Pl(\Sigma^\alpha) > 0$, it holds $k = \tau \, x_0 + o(\tau)$ for some $x_0 > 0$ and all sufficiently large $\tau \in \mathbb{N}$. Using now \eqref{3}, we finally obtain
	\begin{equation*}\label{eq:PLB_2}
	Pl(\Sigma^\alpha) = 8 x_0 \left(\cos \theta(x_0)\right)^2 = 8 \, \sup_{x>0} \left\{ x \left(\cos \theta(x)\right)^2 \right\} > 0.
	\end{equation*}	
	Thus, we see that $Pl(\Sigma^\alpha) = Pl(B)$ for any admissible $\alpha \in (0,2\pi]$. However, the value of $\lim\limits_{\nu \to \infty}\frac{k}{\nu}$ depends on $\alpha$ as follows:
	$$
	\lim_{\nu \to \infty}\frac{k}{\nu} 
	= 
	\frac{\pi}{\alpha} \lim_{\nu \to \infty}\frac{k}{\tau}
	= 
	\frac{\pi x_0}{\alpha}
	=
	\frac{1.165561\ldots}{\alpha}.
	$$

	Finally, let us note that if $\alpha = \pi/m$ for $m \in \mathbb{N}$, then any eigenvalue of \eqref{eq:D} on $\Sigma^\alpha$ has the multiplicity $1$. 
	This fact trivially follows from the validity of Bourget's hypothesis \cite[p.~484]{watson}, since $J_{\nu \frac{\pi}{\alpha}} = J_{\nu m}$, that is, we have the Bessel function of an integer order for any $\nu \in \mathbb{N}$. 
\qed

\subsection{Proof of Proposition \ref{thm:annulus}}

	Consider the ring $A_{r} = \{x \in \mathbb{R}^2:~ r < |x| < 1\}$.
	Separating the variables, it is not hard to find a set of eigenfunctions (up to rotation) of \eqref{eq:D} on $A_r$ in the form
	\begin{equation}\label{eq:eigenfunction}
	\psi_{\nu,k}(\varrho,\vartheta) = \left(J_\nu(a_{\nu,k} \varrho) Y_\nu(a_{\nu,k}) - J_\nu(a_{\nu,k}) Y_\nu(a_{\nu,k} \varrho) \right)\cos(\nu \vartheta),
	\quad 
	\nu \in \mathbb{N} \cup \{0\},~
	k \in \mathbb{N},
	\end{equation}
	and $\lambda_{\nu,k} = a_{\nu,k}^2$ is the eigenvalue associated with $\psi_{\nu,k}$. Here $a_{\nu,k}$ is the $k$-th positive zero of the following cross-product of Bessel functions of the first and second kind:
	\begin{equation}\label{eq:crossprod}
	J_\nu(r z) Y_\nu(z) - J_\nu(z) Y_\nu(r z).
	\end{equation}
	That is,
	$$
	J_\nu(r a_{\nu,k}) Y_\nu(a_{\nu,k}) - J_\nu(a_{\nu,k}) Y_\nu(r a_{\nu,k}) = 0. 
	$$
	It is possible to show that the set of $\psi_{\nu,k}$'s forms an orthogonal basis of $L^2(A_r)$.
	Any eigenvalue $\lambda_{0,k}$ has the multiplicity \textit{at least} $1$ and any other eigenvalue has the multiplicity \textit{at least} $2$. 
	Moreover, $\mu(\psi_{0,k}) = k$ and $\mu(\psi_{\nu,k}) = 2 \nu k$ for $\nu \in \mathbb{N}$. 
	As in the proof of Theorem \ref{thm:disk}, we will write  $\mu(\psi_{\nu,k}) = (2 \nu + \sigma(\nu))k$, where  $\sigma(0)=1$ and $\sigma(\nu)=0$ for $\nu \in \mathbb{N}$.
	
	Assume that the multiplicity of any eigenvalue is at most $2$. As a consequence, we see that \text{any} eigenfunction of \eqref{eq:D} on $A_r$ is of the form \eqref{eq:eigenfunction}, up to scaling. 
	Putting $\lambda_{\nu,k}$ in the increasing order as $\lambda_1 \leq \dots \lambda_n \leq \dots$, noting that $|A_{r}| = \pi (1-r^2)$, and applying the Weyl law, we get
	$$
	n = \lambda_n \frac{|B| |A_{r}|}{4 \pi^2} + o(\lambda_n) = \lambda_n \frac{1-r^2}{4} + o(\lambda_n) = 
	a_{\nu_n,k_n}^2 \frac{1-r^2}{4} + o(a_{\nu_n,k_n}^2).
	$$
	Therefore, passing to a maximizing subsequence and omitting the subindex for $(\nu_n,k_n)$, we deduce that 
	$$
	Pl(A_{r}) = \lim_{\nu+k \to \infty} \frac{4 (2 \nu + \sigma(\nu)) k}{a_{\nu,k}^2 (1-r^2)}.
	$$
	
	Let us note that the inequality \eqref{eq:mccann} from  \cite[Corollary, p.~102]{mccann} used in the proof of Theorem \ref{thm:disk} is, in fact, a consequence of the inequality \cite[(10), p.~102]{mccann} for $a_{\nu,k}$:
	$$
	\left(\frac{a_{\nu,k}}{\nu}\right)^2 
	\geq 
	\left(\frac{a_{0,k}}{\nu}\right)^2 + 1.
	$$
	If, for the maximizing sequence of $Pl(A_r)$, $k$ is bounded, then we can find $C>0$ such that $a_{0,k} > Ck$ for all $k$. Otherwise, if $k \to \infty$, then we use the following  approximation of $a_{0,k}$ of McMahon \cite[(24), p.~29]{mcmahon} (see also Cochran \cite[Theorem, p.~583]{cochran}):
	\begin{equation}\label{eq:a0>a}
	a_{0,k} = \frac{\pi k}{1-r} + O\left(\frac{1}{k}\right).
	\end{equation}
	Therefore, in both cases, we have the existence of $C>0$ such that
	\begin{equation}\label{eq:a>a}
	a_{\nu,k}^2 \geq C k^2 + \nu^2
	\quad 
	\text{for all admissible }
	\nu, k \in \mathbb{N}.
	\end{equation}
	
	Assume first that $Pl(A_r)>0$. Then the case $\nu=0$ can be discarded via \eqref{eq:a0>a}, which implies that $\nu > 0$ and $\mu(\psi_{\nu,k})=2\nu k$.
	Estimating \eqref{eq:a>a} as $a_{\nu,k}^2 \geq C k^2$ and as $a_{\nu,k}^2 \geq \nu^2$, we obtain
	$$
	Pl(A_r) \leq \frac{8}{1-r^2} \lim_{\nu + k \to \infty} \min\left\{
	\frac{k}{\nu}, \frac{\nu}{Ck}
	\right\}.
	$$
	Since $Pl(A_r)>0$, we conclude that both $\nu$ and $k$ tend to infinity, and there exist $A_1, A_2 > 0$ such that 
	$$
	A_1 k < \nu < A_2 k
	\quad 
	\text{for all sufficiently large }
	k \in \mathbb{N}.
	$$
	Therefore, up to a subsequence, we have $\nu = k \, x_0 + o(k)$ for some $x_0>0$ and all large $k \in \mathbb{N}$. 
	Thus, 	
	\begin{equation}\label{eq:upper_bound_pa}
	Pl(A_{r}) = \frac{8}{1-r^2} \, \sup_{x>0} \left\{x \,\limsup_{k \to \infty} \frac{k^2}{a_{kx,k}^2}\right\}.
	\end{equation}
	
	Assume now that $Pl(A_r)=0$. Then the expression on the right-hand side of \eqref{eq:upper_bound_pa} provides a lower bound for $Pl(A_r)$, which again implies that the equality \eqref{eq:upper_bound_pa} is satisfied.
\qed

\begin{cor}
	Combining the upper estimate \eqref{eq:Pleijel} (which is strict by \cite{B,stein}) with the expression \eqref{eq:upper_bound_pa}, we obtain the following lower bound for the $k$-th positive zero of the cross-product of Bessel functions \eqref{eq:crossprod} of order $kx$, $x>0$:
	$$
	a_{kx,k} > \frac{\sqrt{2} j_{0,1} k}{\sqrt{1-r^2}} + o(k) > 
	\frac{3.4 \,k}{\sqrt{1-r^2}} + o(k)
	$$
\end{cor}

\subsection{Proof of Proposition \ref{thm:annulus_sector}}

	Consider the annular sector 
	$$
	\Sigma^\alpha_r = 
	\left\{(\varrho \cos\vartheta, \varrho \sin\vartheta):~ \varrho \in (r,1),~ \vartheta \in \left(-\frac{\alpha}{2},\frac{\alpha}{2}\right)
	\right\},
	\quad
	r \in (0,1), ~
	\alpha \in (0,2\pi].
	$$
	By separation of variables, it can be derived that there is a basis of eigenfunctions of \eqref{eq:D} on $\Sigma^\alpha_r$ of the form
	$$
	\psi_{\nu,k}^\alpha(\varrho,\vartheta) 
	= \left(J_{\nu\frac{\pi}{\alpha}}(a_{\nu\frac{\pi}{\alpha},k} \varrho) Y_{\nu\frac{\pi}{\alpha}}(a_{\nu\frac{\pi}{\alpha},k}) 
	- 
	J_{\nu\frac{\pi}{\alpha}}(a_{\nu\frac{\pi}{\alpha},k}) Y_{\nu\frac{\pi}{\alpha}}(a_{\nu\frac{\pi}{\alpha},k} \varrho)\right)
	\sin\left(\nu \pi \left(\frac{\vartheta}{\alpha}+\frac{1}{2}\right)\right)
	$$
	for $\nu, k \in \mathbb{N}$, 
	with the associated eigenvalues $\lambda_{\nu,k} = a_{\nu\frac{\pi}{\alpha},k}^2$. Here $a_{\nu\frac{\pi}{\alpha},k}$ is the $k$-th positive zero of the cross-product of Bessel functions \eqref{eq:crossprod} of order $\nu\frac{\pi}{\alpha}$.
	It is not hard to see that $\mu(\psi_{\nu,k}^\alpha) = \nu k$, and any eigenvalue $\lambda_{\nu,k}$ has the multiplicity at least $1$.
	
	Assuming that each $\lambda_{\nu,k}$ has the multiplicity exactly $1$, we	argue in much the same way as in the proof of Theorem \ref{thm:annulus} (and Theorem \ref{thm:sector-disk}) to conclude that
	\[
	Pl(\Sigma^\alpha_{r}) = 
	Pl(A_r)
	=
	\frac{8}{1-r^2} \, \sup_{x>0} \left\{x \,\limsup_{k \to \infty} \frac{k^2}{a_{kx,k}^2}\right\}.
	\pushQED{\qed}
	\qedhere
	\]

\subsection{Proof of Lemma \ref{lem:multiplicity}}
\begin{figure}[ht]
	\centering
	\begin{minipage}[t]{0.24\linewidth}
		\centering
		\includegraphics[width=0.99\linewidth]{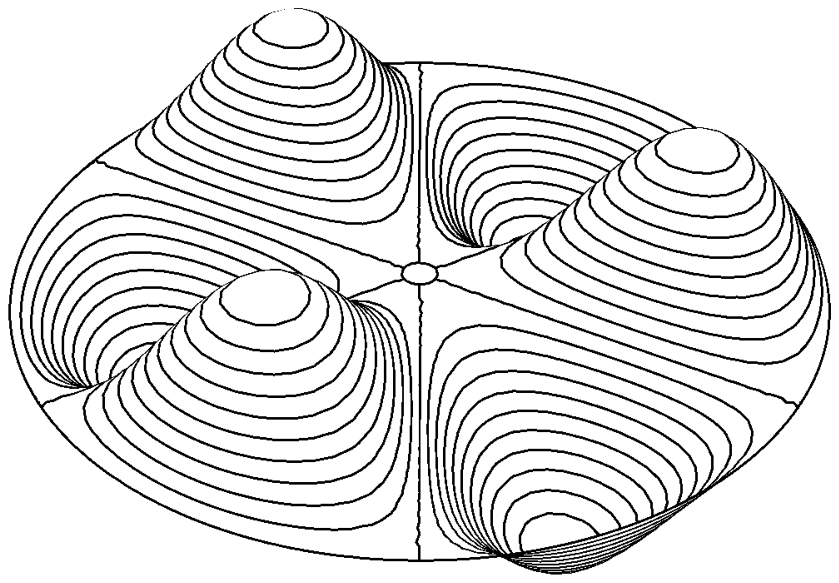}\\
		$C_1=1$, $C_2=0$
	\end{minipage}
	\begin{minipage}[t]{0.24\linewidth}
		\centering
		\includegraphics[width=0.99\linewidth]{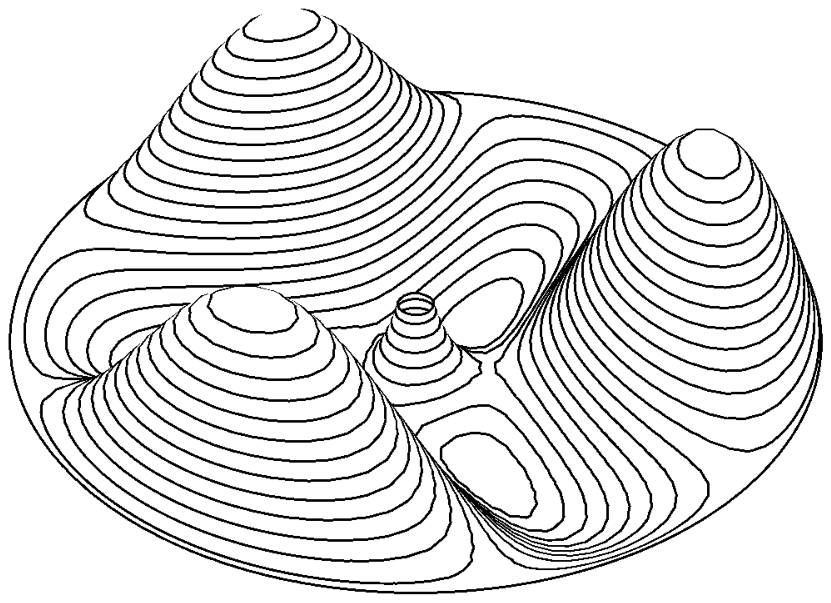}\\
		$C_1=1$, $C_2=0.5$
	\end{minipage}
	\begin{minipage}[t]{0.24\linewidth}
		\centering
		\includegraphics[width=0.99\linewidth]{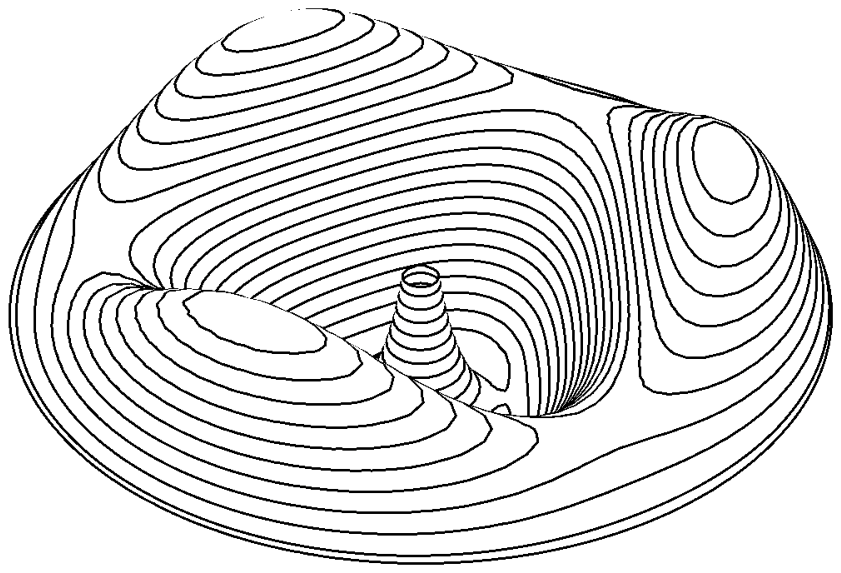}\\
		$C_1=0.5$, $C_2=1$
	\end{minipage}
	\begin{minipage}[t]{0.24\linewidth}
		\centering
		\includegraphics[width=0.99\linewidth]{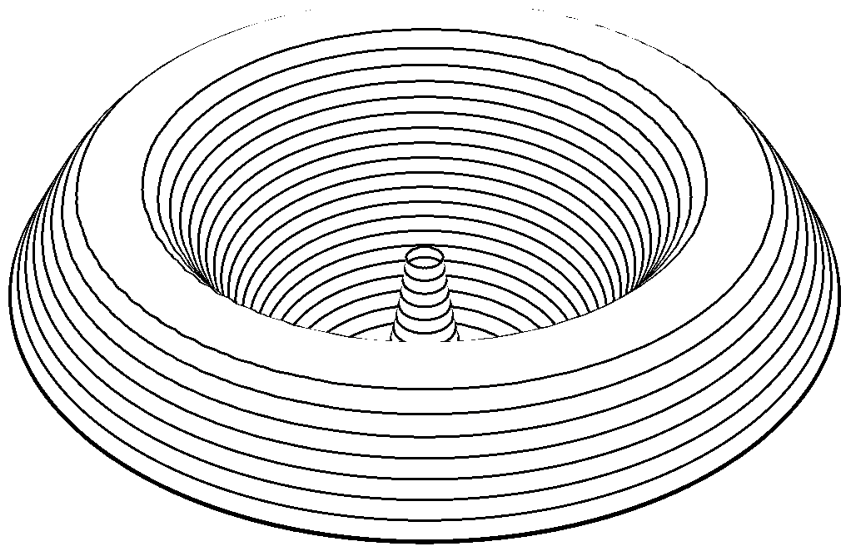}\\
		$C_1=0$, $C_2=1$
	\end{minipage}
	\caption{Eigenfunctions $C_1\psi_{3,1} + C_2\psi_{0,2}$ of \eqref{eq:D} on $A_{r_0}$ with $r_0 \approx 0.044951$ associated with $\lambda = \lambda_{3,1}=\lambda_{0,2}$.}
	\label{fig1}
\end{figure}

	Note that any zero $a_{\nu,k}(r)$ of the cross-product of Bessel functions \eqref{eq:crossprod} of order $\nu$ is a continuous (in fact, analytic) function with respect to $r > 0$, see Cochran \cite[Theorem 1]{cochran2}. Here $r>0$ stands for the inner radius of the ring $A_r$.
	Then, the numerical computations show that $a_{0,2}(0.01) \approx 6.0109$ and $a_{3,1}(0.01) \approx 6.3801$, while $a_{0,2}(0.1) \approx 6.8575$ and $a_{3,1}(0.1) \approx 6.3804$.
	Therefore, by continuity, there exists $r_0 \in (0.01,0.1)$ such that $a_{3,1}(r_0) = a_{0,2}(r_0)$, which implies that the eigenvalue $\lambda = \lambda_{3,1} = \lambda_{0,2}$ has the multiplicity at least $3$. (According to the numerical computations, $r_0 \approx 0.044951$ and $\lambda \approx 40.7064$; see Fig.~\ref{fig1}.)
	Actually, this fact was observed already by Kline \cite[Fig.~1]{klin}.
	The claim for annular sectors can be obtained either by the same method as above or by the arguments of \cite[Section 2.3]{BL}.
\qed

\bigskip
\noindent
{\bf Acknowledgements.}
The author wishes to thank Bernard Helffer for a stimulating discussion and helpful remarks. 
This research has been supported by the Grant Agency of the Czech Republic, project 18-03253S, and by the project LO1506 of the Czech Ministry of Education, Youth and Sports.

\end{document}